\documentclass[a4paper,11pt]{article}
\usepackage{graphicx, amsmath, amsthm, amssymb, mathrsfs}
\usepackage{a4wide, enumerate}
\usepackage[utf8]{inputenc}
\usepackage{slashbox}
\DeclareMathOperator{\erf}{erf}

\newcommand{\mbR}{\mathbb{R}}  % Real numbers
\newtheorem{definition}{Definition}[section]
\newtheorem{remark}{Remark}[section]
\newtheorem{theorem}{Theorem}[section]

\newtheorem{problem}{Problem}[section]
\allowdisplaybreaks
\title{Parameter Estimation of Fiber Lay--down in Nonwoven Production \\-- An Occupation Time Approach --}
\author{Wolfgang Bock$^{\dagger}$,  Thomas Götz$^{\star}$, Uditha Prabhath Liyanage$^{\dagger}$\\[0.5cm]
 $^{\dagger}$ \small{Dept.~of Mathematics, University of Kaiserslautern,}\\
\small{P.O.Box ~3049, D--67653 Kaiserslautern, Germany.}\\
\small{\emph{$\{$bock, liyanage$\}$@mathematik.uni-kl.de}}\\[0.5 cm]
 $^{\star}$ \small{Mathematisches Institut, Universit\"at Koblenz, Universit\"atsstr. 1}\\
\small{D-56070 Koblenz, Germany.}\\
\small{\emph{goetz@uni-koblenz.de}}
}

\begin{document}
\maketitle
\begin{abstract}
 {In this paper we investigate the parameter estimation of the fiber lay--down process in the production of nonwovens. The parameter estimation is based on the mass per unit area data, which is available at least on an industrial scale. We introduce a stochastic model to represent the fiber lay--down and through the model's parameters we characterize this fiber lay--down. Based on the occupation time, which is the equivalent quantity for the mass per unit area in the context of stochastic dynamical systems, an optimization procedure is formulated that estimates the parameters of the model. The optimization procedure is tested using occupation time data given by Monte--Carlo simulations. The feasibility of the optimization procedure on an industrial level is tested using the fiber paths simulated by the industrial software FYDIST. }
\end{abstract}

\section{Introduction}
\begin{figure}[h]
\begin{center}
\includegraphics[width=.45\textwidth]{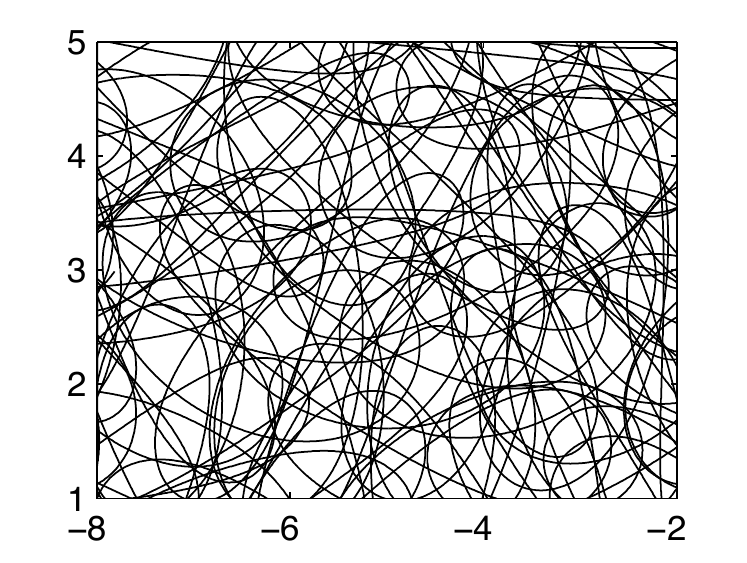}\hfill
\includegraphics[width=.45\textwidth]{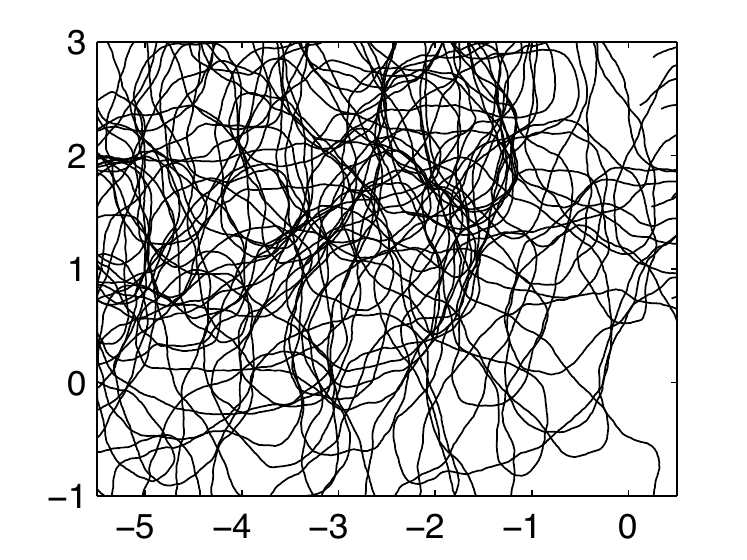}\\
\caption{Forms of fiber webs (nonwovens).}\label{fig:webForms}
\end{center}
\end{figure}

Nonwoven materials or fleece are webs of long flexible fibers that are used for composite materials, e.g.~filters, as well as in the hygiene and textile industries.  They are produced in melt--spinning operations: hundreds of individual, endless fibers are obtained by the continuous extrusion of a molten polymer through narrow nozzles, which are densely and equidistantly placed in a row at a spinning beam. The viscous or viscoelastic fibers are stretched and spun until they solidify due to cooling air streams. Before the elastic fibers lay down on a moving conveyor belt to form a web, they become entangled and form loops due to highly turbulent air flows. Figure \ref{fig:webForms} shows  at a microscopic level the webs formed on the conveyor belt.
The homogeneity and load capacity of the fiber web are the most important textile properties for quality assessment of industrial nonwoven fabrics. The optimization and control of the fleece quality require modeling and simulation of the fiber dynamics and the lay--down.

There are two classes of approaches to model the fiber lay--down process. The first class uses microscopic details to model the lay--down of the fibers. Here, each fiber is seen as an elastic beam. The software FYDIST, developed by the Fraunhofer ITWM, Kaiserslautern, Germany, uses such models to describe the fiber lay--down. Since the motion of each fiber is simulated using the physics on a microscopic level, the behavior of the simulated fibers is quite close to the real industrial fibers. Nevertheless, due to the large number of microscopic details included, the numerical computations are highly time consuming. 

The second approach is based on macroscopic, quantitative description. Here, the lay--down is modeled stochastically, i.e. the models are consist of stochastic differential equations with a certain set of parameters. Since this quantitative approach does not use fine details of the fiber lay--down, it allows for fast numerical simulations.

\begin{figure}[h]
\begin{center}
\includegraphics[width =0.45\textwidth]{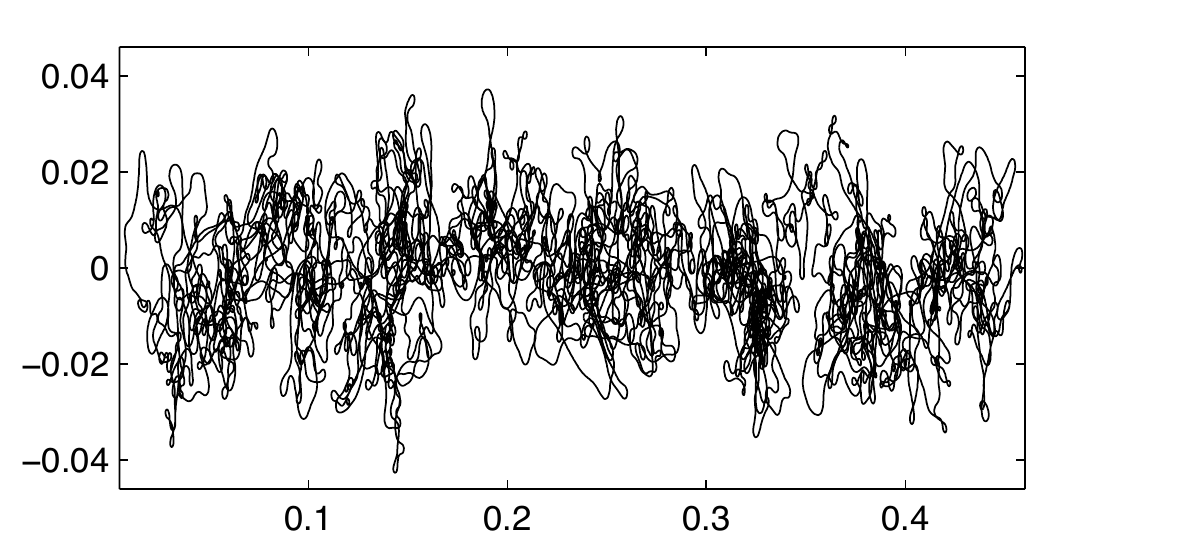}\hfill
\includegraphics[width = 0.45\textwidth]{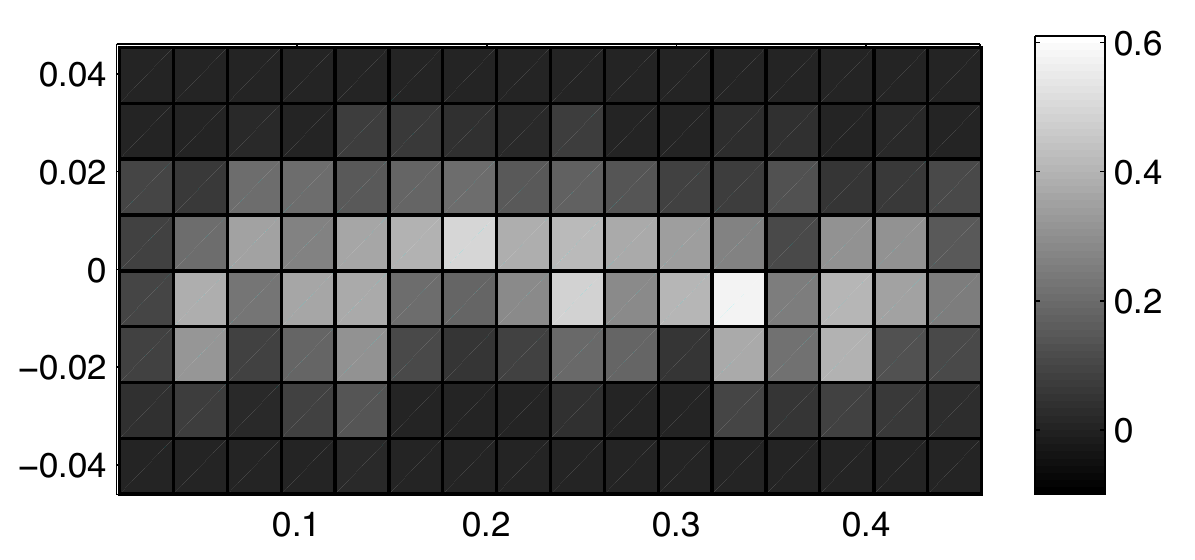}\\
\caption{A fiber path (left) laid on a conveyor belt and the mass distribution of the fiber (right). }\label{fig:MassDistribution}
\end{center}
\end{figure}

Available data to judge the quality, at least on the industrial scale, are usually the mass per unit area of the fleece. Figure \ref{fig:MassDistribution} (right) shows the mass distribution per area of a fleece made by a single fiber (left).
Since the mass distribution is an averaged quantity, the microscopic details of the fibers do not play a significant role. Nevertheless, the quantitative, macroscopic approach provides sufficient details to optimize fibers based on mass distributions. Moreover, due to fast simulations the stochastic models increase the efficiency of the  optimization. Thus, in order to identify the parameter of the fiber lay--down we use the macroscopic approach.\\

A stochastic model for the fiber deposition in the nonwoven production was proposed and analyzed in Ref.~\cite{BGKMW08,GKMW07}. In particular, the stochastic model proposed in Ref.~\cite{BGKMW08}, which represents the fiber deposition on a moving belt, is constituted by nonlinear stochastic differential equations (SDE). \\
The hydrodynamic scaling limit of the resulting stochastic process of the model is given by an Ornstein--Uhlenbeck process with moving mean. The aim of this paper is to determine the parameters of the Ornstein--Uhlenbeck process with moving mean from available mass per unit area data, i.e.~the occupation time.

The paper is organized as follows: In Section~\ref{sec:ModelandTheory} we introduce the fiber lay--down process as a simplified model for the fiber deposition. Furthermore, we define the occupation time and present an analytical expression of the expected occupation time. An optimization procedure to estimate the model parameters from available occupation times is presented along with numerical experiments in Section~\ref{sec:Numerics}. Finally, we draw some conclusions and give an outlook to open questions.

\section{Model and theory}\label{sec:ModelandTheory}
Summarizing \cite{BGKMW08}, we model the fiber lay--down process on a moving conveyor belt by a stochastic process $Y:=(Y_t^1,Y_t^2)^T_{t\geq 0}\in\mbR^2$ satisfying the stochastic differential equations 
\begin{equation}
 \begin{split}
 d Y_t^1 &= \lambda(\theta_1+\kappa t-Y_t^1)\, dt + \sigma_1 \, d W_t^1,\\[0.5ex]
 d Y_t^2 &= \lambda(\theta_2-Y_t^2)\, dt + \sigma_2 \, d W_t^2,\quad\quad 
 \end{split}
 \end{equation}
with initial condition $Y_0 = (Y^1_0,Y^2_0)^T\in\mbR^2$.  The stiffness $\lambda>0$ governs the deterministic part of the fiber deposition, while a standard two--dimensional Brownian motion $W=(W_t^1,W_t^2)^T_{t\geq 0}$ and the diffusion parameters $\sigma_1,\,\sigma_2\in\mbR$ govern the stochastic part. The parameter $\kappa\ge 0$ is the belt speed and $(\theta_1,\theta_2)^T\in\mbR^2$ is the reference position of the fiber lay--down.\\
The fiber lay--down process is always centered to its moving position $(\theta_1+\kappa t,\theta_2)$. Assuming the deterministic forces are symmetric around the reference position of the fiber lay--down, we use the same stiffness coefficient $\lambda$ in both equations, i.e. an isotropic fiber lay--down. For sake of simplicity, we set $\theta_1=\theta_2=0$. Then, the fiber lay--down model reads as
\begin{equation}\label{eqn:twoDimMovingOUmodel}
 \begin{split}
 d Y_t^1 &= -\lambda(Y_t^1-\kappa t)\, dt + \sigma_1 \,d W_t^1,\\[0.5ex]
 d Y_t^2 &= -\lambda\,Y_t^2\,dt + \sigma_2 \,dW_t^2.
 \end{split}
 \end{equation}
The random variable $Y_t$ models the deposition point of an individual fiber on the fleece.  If we follow the random variable over a time interval $[0,T]$ for $T>0$, we obtain the path of an individual fiber. The fiber lay--down process $Y$ on a moving belt is characterized by the drift and diffusion parameters $\lambda$, $\sigma_1$ and $\sigma_2$, together with the belt speed $\kappa$. Assuming the belt speed $\kappa$ as a known information, the process  $Y$ given by the parameters $\lambda$, $\sigma_1$ and $\sigma_2$ may be denoted as $Y_{\lambda,\sigma_1,\sigma_2}$.

To introduce the mathematical analogue of the mass per unit area we need the following definition.
\begin{definition}[Occupation time]\label{def:occ_time}
Let $T>0$ and consider a rectangle  $\mathcal{D}:=[a_1,b_1]\times[a_2,b_2]\subset\mbR^2$. The \emph{occupation time} $M_{\mathcal{D},\,T}$ is defined as
\begin{equation*}
    M_{\mathcal{D},\,T}(Y) := \int_0^T \mathbf{1}_{\mathcal{D}}(Y_t) \, dt .
 \end{equation*}
Here, $\mathbf{1}_{\mathcal{D}}$ denotes the indicator function of the rectangle $\mathcal{D}$.
\end{definition}

\begin{remark}
The occupation time is a random variable itself. It models the time, the random process spends inside the rectangle $\mathcal{D}$ during the time interval $[0,T]$. In terms of our physical model for the nonwoven production, the occupation time can be interpreted as the mass of fiber material deposited inside $\mathcal{D}$, i.e.~the mass per area of the final fleece. This quantity is easily accessible even on the scale of industrial production and hence it will serve as the input to our parameter estimation problem.
\end{remark}

The following theorem can be proven using the techniques of white noise analysis, see \cite{BGGL2011} for the one--dimensional case.
\begin{theorem}[Occupation time]\label{thm:expectedOccMovingBelt}
Let $M_{[a_1,b_1]\times[a_2,b_2],\,T}(Y)$ be the occupation time of the fiber lay--down process $Y$ given by \eqref{eqn:twoDimMovingOUmodel}. Then, its expectation is given by 
\begin{multline}\label{eqn:ExpOccuTimeMoving}
\mathbb{E}(M_{[a_1,b_1]\times [a_2,b_2],\,T}(Y))\\
= \frac{1}{4}\int\limits_0^T\left({\erf}\left(\frac{\sqrt{\lambda}\left(b_1-\frac{\kappa}{\lambda}\left( \lambda\,t-1+e^{-\lambda\,t} \right)\right)}{\sigma_1\sqrt{1-e^{-2 \lambda t}}}\right)\right.- 
\left.{\erf}\left(\frac{\sqrt{\lambda}\left(a_1-\frac{\kappa}{\lambda}\left( \lambda\,t-1+e^{-\lambda\,t} \right)\right)}{\sigma_1\sqrt{1-e^{-2 \lambda t}}}\right)\right) \\
  \times\left({\erf}\left(\frac{\sqrt{\lambda}b_2}{\sigma_2\sqrt{1-e^{-2 \lambda t}}}\right)-{\erf}\left(\frac{\sqrt{\lambda}a_2}{\sigma_2\sqrt{1-e^{-2 \lambda t}}}\right)\right) \, dt.
\end{multline}
\end{theorem}
In the following we use $E_{[a_1,b_1]\times [a_2,b_2],\,\kappa,\,T}\-(\lambda,\sigma_1,\sigma_2)$ to denote the expected occupation time of the process $Y_{\lambda,\sigma_1,\sigma_2}$.

\section{Numerics}\label{sec:Numerics}

\subsection{Monte--Carlo computation of the expected occupation time}
\begin{figure}[h]
\begin{center}
\includegraphics[width =0.45\textwidth]{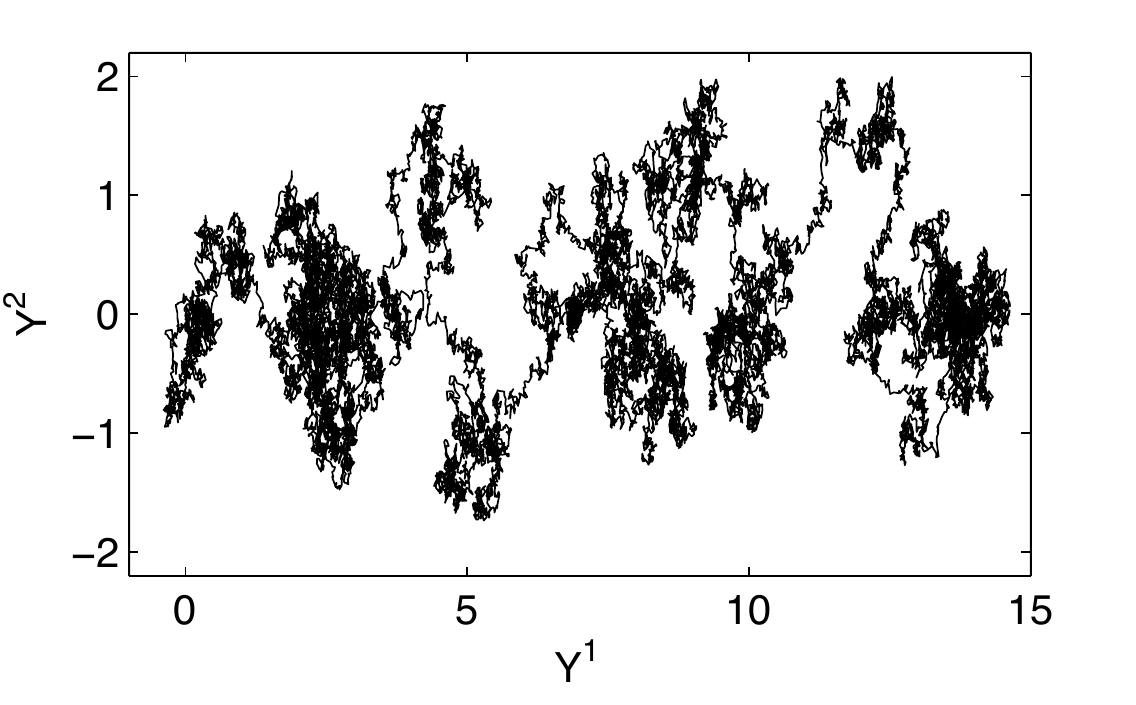}
\includegraphics[width =0.45\textwidth]{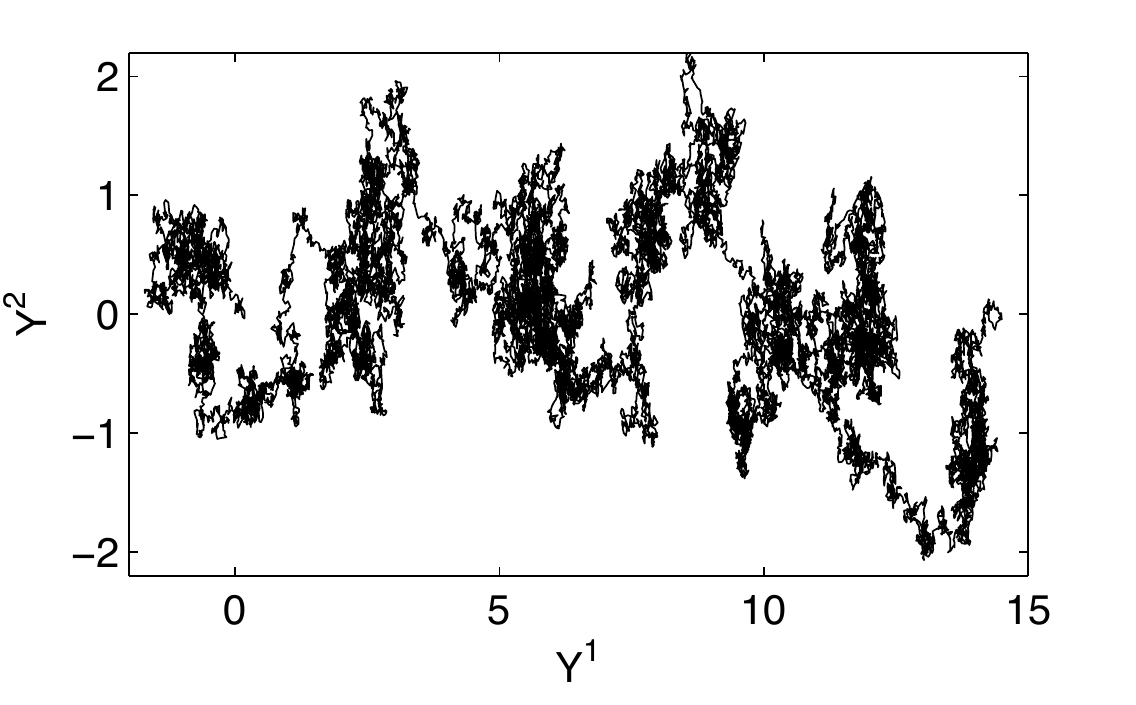}
\caption{Sample path of the process $Y$ given by the parameters $\lambda=1$, $\sigma_1=1$ and $\sigma_2 = 1$. The time horizon and the belt speed are  $T=30$ and $\kappa =0.5$, respectively.}\label{fig:movingSamplePath}
\end{center}
\end{figure}

A Monte--Carlo method that approximates the expected occupation time of the fiber lay--down process $Y$ consists of the following steps
\begin{itemize}
\item approximate numerically the process $Y$ given by the SDE \eqref{eqn:twoDimMovingOUmodel},
\item approximate  the expected occupation time by considering a sufficiently large number of sample paths of a two--dimensional Brownian motion $W$.
\end{itemize}
To approximate the process numerically we use the standard Euler--Maruyama method up to time horizon $[0,T]$, see \cite{kloeden92}. Let us denote the Euler--Maruyama approximation on the time grid $\{0=\tau_0<\tau_1<\dots<\tau_k=T\}$ by $Y^{(k)}$. Figure \ref{fig:movingSamplePath} shows two sample paths of the process $Y$. Both sample paths are plotted with respect to the Cartesian coordinate system defined by the nozzle position and the moving direction of the conveyor belt. The axis that passes through the nozzle position and parallel to the moving direction of the conveyor belt is called the central axis of the fiber lay--down.

Let $\tilde{\Phi}$ denote the mapping implicitly defined by the stochastic differential equations~\eqref{eqn:twoDimMovingOUmodel}, i.e. $Y=\tilde{\Phi}(W)$ for a given sample path of the Brownian motion $W$. Let $\tilde{\Phi}^{(k)}$ be the mapping implicitly given by the Euler--Maruyama scheme, i.e. $Y^{(k)}=\tilde{\Phi}^{(k)}(W)$. Thus, computing the integral
\begin{equation*}
M_{{[a_1,b_1]\times [a_2,b_2]},\,\kappa,T}(\tilde{\Phi}^{(k)}(W)) =\int_0^T\textbf{1}_{[a_1,b_1]\times [a_2,b_2]}(\tilde{\Phi}^{(k)}(W))\, dt,
\end{equation*}
we get the occupation time of $Y^{(k)}$. The Monte--Carlo approximation for the expected occupation time ${E}\-(M_{{[a_1,b_1]\times [a_2,b_2]}\-,\,\kappa,\,T}\-(Y)) $ is given by
\begin{equation}\label{eqn:movingMCApproximation}
S_N(\tilde{\Phi}^{(k)}(W)) = \frac{1}{N}\sum_{j=1}^N M_{{[a_1,b_2]\times[a_2,b_2]},\,\kappa\,,T}(\tilde{\Phi}^{(k)}(W_j)),
\end{equation}
where $W_1,W_2,\-...,W_N$ are independent sample paths of the two--dimensional Brownian motion $W$. Due to the \emph{law of large numbers} we know that $S_N(Y^{(k)})\-  \rightarrow\-{E}\-(M_{{[a_1,b_1]\-\times [a_2,b_2]},\,\kappa,\,T}(Y))$ as $N\rightarrow \infty$. For more details see \cite{MuellerRitter2009} and references therein.

Table \ref{tab:movingTable01} lists computations of the expected occupation times for the process $Y$ corresponding to the area $[3,15]\times[-1,1]$ for different parameters $\lambda$, $\sigma_1$, $\sigma_2$, belt speeds $\kappa$ and time horizons $[0,T]$. The Monte--Carlo method uses $10000$ independent sample paths.
%\marginpar{Why are the values for $T=50$ smaller than for $T=30$??? This cannot be, since occupation increases monotonically with $T$. Please check!}
\begin{table} [h]
\begin{center}
\begin{tabular}{|c|c|c||c|c|}                                                                                                                                                                                                                                                                                                                                                                                                                                                                                                                                                                                                                                                                                                                                                                                                                                                                                                                                                                                                                                                                                                                                                                                                                                                                                                                                                                                                                                                                                                                                                                                                                                                                                                                                                                                                                                                                                                                                                                                                                                                                                                                                                                                                                                                                                                                                                                                                                                                                                                                                                                                                                                                                                                                                                                                                                                                                                                                                                                                                                                                                                                                                                                                                                                                                                                                                                                                                                                                                                                                                                                                                                                                                                                                                                                                                                                                                                                                                                                                                                                                                                                                                                                                                                                                                                                                                                                                                                                                                                                                                                                                                                                                                                                                                                                                                                                                                                                                                                                                                                                                                                                                                                                                                                                                                                                                                                                                                                                                                                                                                                                                                                                                                                                                                                                                                                                                                                                                                                                                                                                                                                                                                                                                                                                                                                                                                                                                                                                                                                                                                                                                                                                                                                                                                                                                                                                                                                                                                                                                                                                                                            \hline
 & $\kappa=1$    &$\kappa=2$    &$\kappa=1$   & $\kappa=2$    \\
  				\backslashbox{Time $T$}{Parameters}		 	&	$\lambda,\sigma_1,\sigma_2=1$  	&     $\lambda,\sigma_1,\sigma_2=1$ &	$\lambda,\sigma_1,\sigma_2=2$      			&       $\lambda,\sigma_1,\sigma_2=2$  \\ \hline
$T = 7$ by MC 													& 			2.5465 		 &		3.8793		& 		2.3899	 &		3.4246	    \\ \hline
$T = 7$ by \eqref{eqn:ExpOccuTimeMoving} 	& 			2.5484  	 &		3.8742		& 		2.3910	 &		3.4207 	    \\ \hline
$T = 30$ by MC 												& 			10.1382 	 &		5.1327  	& 		8.1910  &		4.1021	     \\ \hline
$T = 30$ by \eqref{eqn:ExpOccuTimeMoving} 	& 			10.1327	 &		5.1380  	& 		8.1938	 &  	4.1063	    \\ \hline
$T = 50$ by MC 												& 			10.1382	 &		5.1327		& 		8.1910	 &		4.1021	     \\ \hline
$T = 50$ by \eqref{eqn:ExpOccuTimeMoving} 	& 			10.1327	 &		5.1380		& 		8.1938	 &		4.1063	    \\ \hline
\end{tabular}
\end{center}
\caption{Expected occupation times of $Y$ computed by the Monte--Carlo (MC) method and the formula \eqref{eqn:ExpOccuTimeMoving} .}\label{tab:movingTable01}
\end{table}

Note that for given parameters $\lambda$, $\sigma_1$ and $\sigma_2$ the expected occupation times corresponding to the time horizon $T=30$ do not differ from those corresponding to the time horizon $T=50$. This observation is due to the movement of the belt in the lay--down process. Note that in cases $T=30$ and $T=50$ the Monte--Carlo method uses the fiber sample path corresponding to random sequences generated by a fixed set of seeds.
%we have used the same set of seeds to generate the random sequence that generate the sample paths used in the Monte--Carlo method in the both cases $T=30$ and $T=50$.
%the in cases $T=30$ and $T=50$ the Monte--Carlo method uses the fiber sample path given by random sequences generated by fixed set of seeds. 
%  Moreover, increasing the belt speed, while the other settings are fixed, leads to lower expected occupation times.
 
\begin{figure}[h]
\begin{center}
\includegraphics[width =0.45\textwidth]{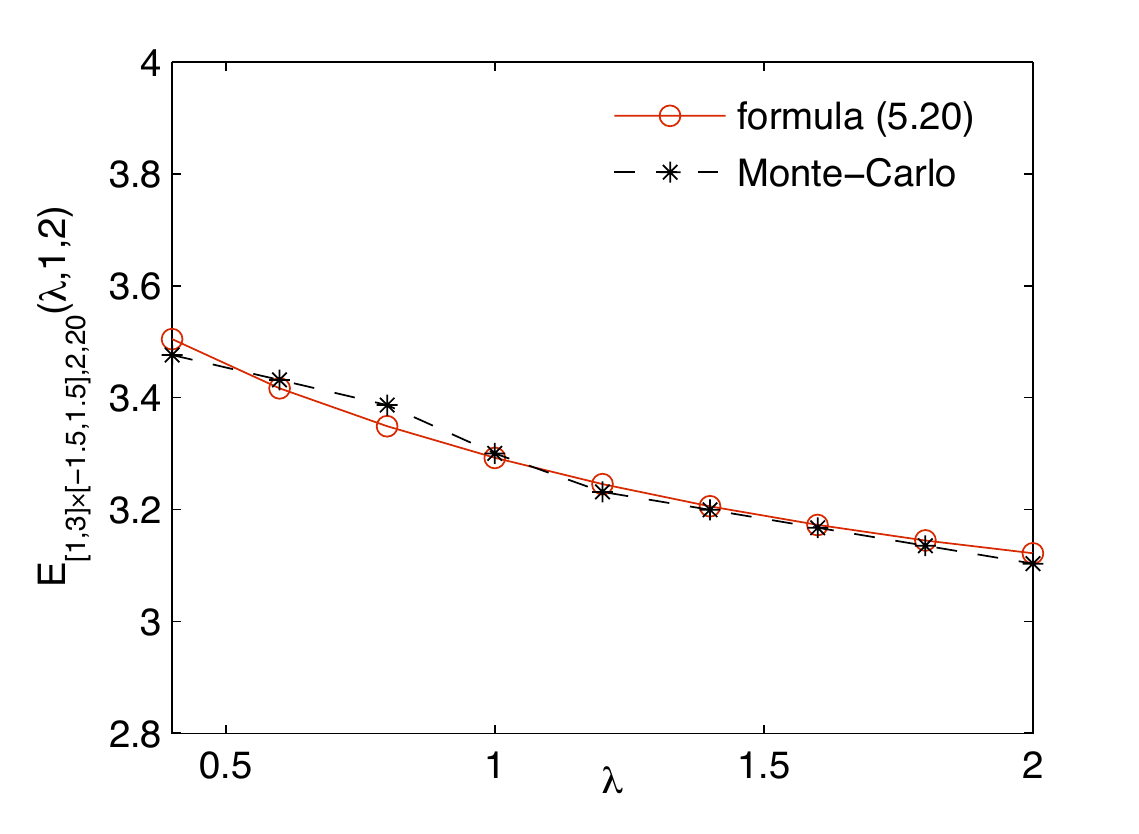} \hfill
\includegraphics[width =0.45\textwidth]{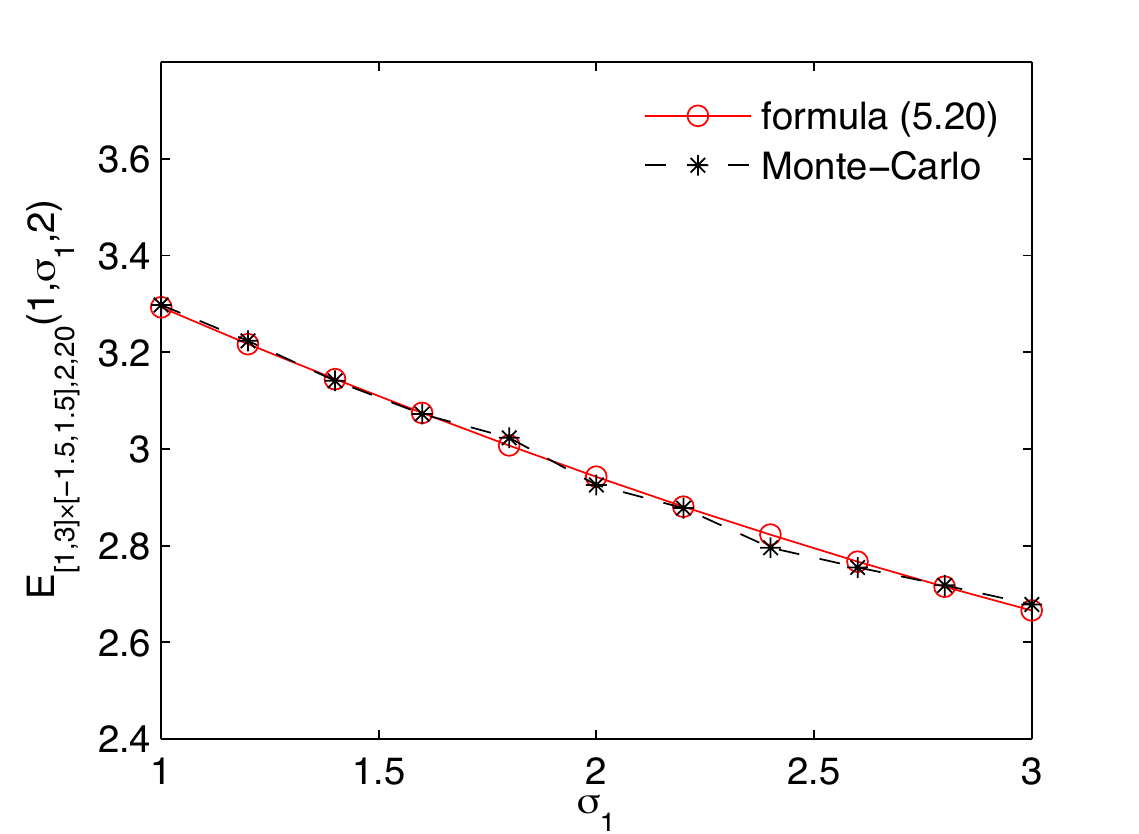}
\caption{Comparison of the expected occupation times $E_{\mathcal{D},\,\kappa , \,T}\-(\lambda,\sigma_1,\sigma_2)$ computed by the formula \eqref{eqn:ExpOccuTimeMoving} and the Monte--Carlo method for the domain $\mathcal{D}=[1,3]\times[-1.5,1-5]$, belt speed $\kappa=2$ and time horizon $T=20$.}\label{fig:ComparisonexpectedOccupationTimeMoving}
\end{center}
\end{figure}

Figure \ref{fig:ComparisonexpectedOccupationTimeMoving} compares the computation of expected occupation times computed by the analytical formula \eqref{eqn:ExpOccuTimeMoving} and the Monte--Carlo method with $10000$ sample paths. 

To analyze how the expected occupation time of a given process $Y$ is distributed around the central axis of the fiber lay--down, we proceed as follows:\\
Let us consider a fiber lay--down process $Y$ given by the parameters $\lambda=1.5$, $\sigma_1=1$ and $\sigma_2=1$. We fix the belt speed $\kappa= 1$, time horizon $T = 20$ and consider an area given by $[5,15]\times[-1.8,1.8]$, which is symmetric around the central axis.
\begin{figure}[h]
\begin{center}
\includegraphics[width =0.45\textwidth]{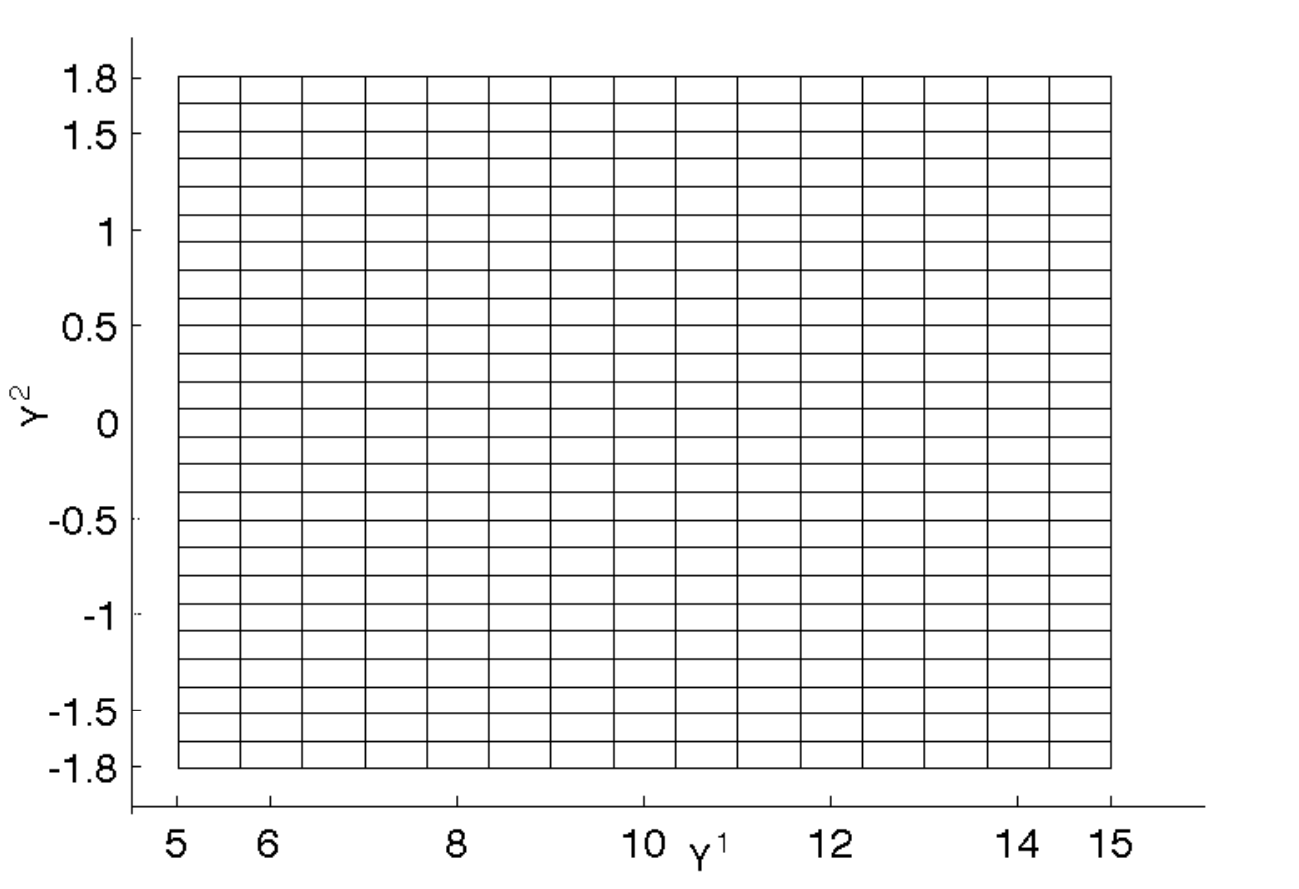} \hfill
\includegraphics[width =0.45\textwidth]{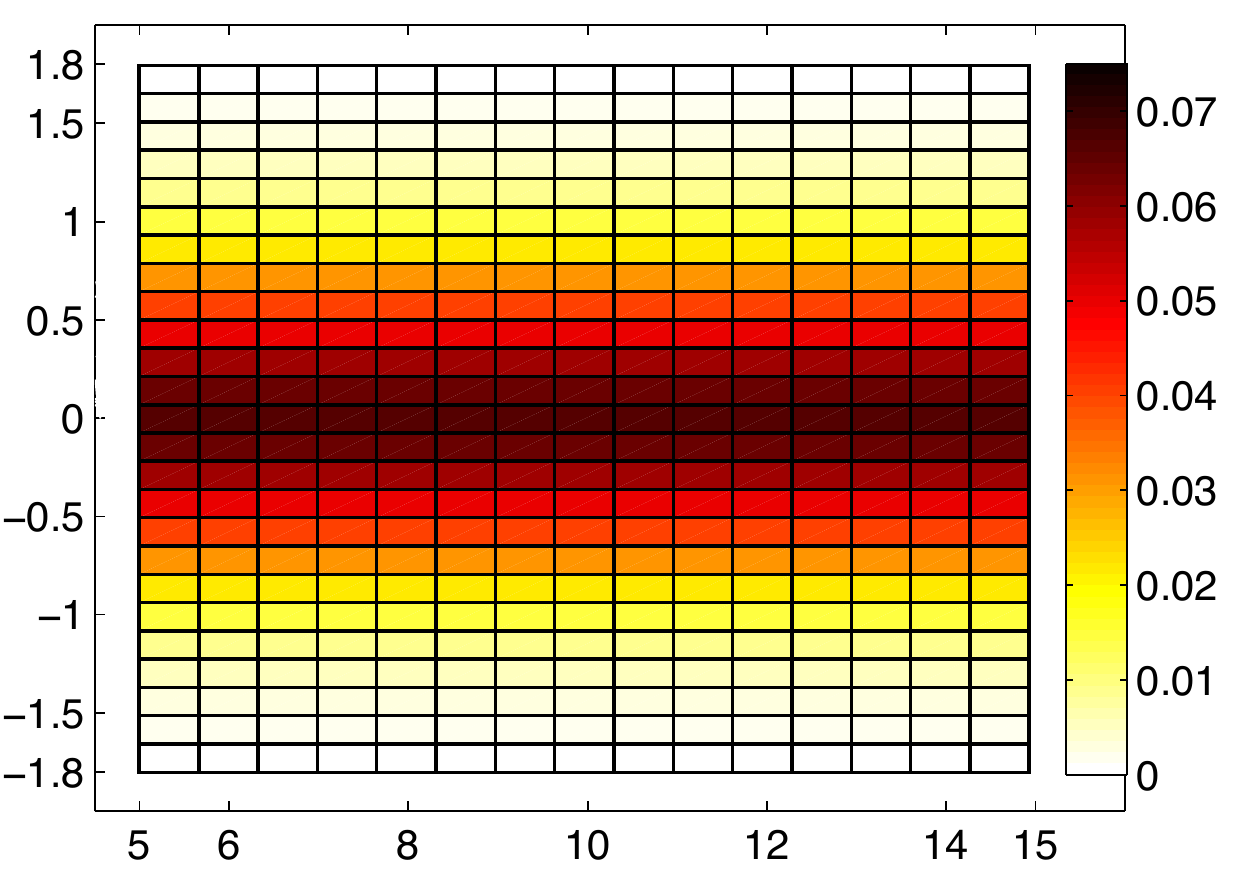}\\
\caption{The grid that divides the area $[5,15]\times[-1.8,1.8]$ (left) and the expected occupation time distribution (right) with respected to the cell areas given by the grid.}\label{fig:expectedOccupationTimeAreawise}
\end{center}
\end{figure}
We subdivide this area into a grid of $25\times15$ rectangular sub areas (cells) of equal sized, see Figure \ref{fig:expectedOccupationTimeAreawise} (left). Next, we compute the expected occupation times corresponding to each cell using the formula \eqref{eqn:ExpOccuTimeMoving}. On the right hand side of Figure \ref{fig:expectedOccupationTimeAreawise} we plot the corresponding expected occupation times. Note that the occupation time decreases as we move away from the central axis $Y^2=0$ of the lay--down process. Increasing the parameter $\sigma_2$ or decreasing the parameter $\lambda$ leads to a flater and broader distribution of the expected occupation time with respect to the $Y^2$--coordinate.

\subsection{An optimization procedure for the parameter estimation}

As we have already mentioned, the fiber lay--down process $Y$ is characterized by the parameters $\lambda$, $\sigma_1$,  $\sigma_2$ and  the belt speed $\kappa$. Using this relation inversely, we estimate the parameters $\lambda$, $\sigma_1$ and $\sigma_2$ based on a given set of expected occupation times. Including the belt speed $\kappa$ as a known information, we formulate the parameter estimation as an optimization problem.

\begin{problem}\label{prob:estimationMovingBelt}
Assume that the occupation times $E_{i,j}$ for $i=1,2,\dots,n$ and $j=1,2,\dots,m$ corresponding to a fixed time horizon $T$, different areas $[a^i_1,b^i_1]\times[a^i_2,b^i_2]$ and various belt speeds $\kappa_j$ are given. Determine the parameters $\lambda,\sigma_1,\sigma_2$ such that the least--squares deviation
\begin{equation}\label{eqn:costFunctionMovingBelt}
 R(\lambda,\sigma_1,\sigma_2):=\sum_{i=1}^n \sum_{j=1}^m\left(E_{[a^i_1,b^i_1]\times[a^i_2,b^i_2],\,\kappa_j ,\, T}(\lambda,\sigma_1,\sigma_2)-E_{i,j}\right)^2,
\end{equation}
is minimal.
\end{problem}

To minimize $R(\lambda,\sigma_1, \sigma_2)$ we use the simplex search method implemented in \textsc{Matlab}  as the function \textit{fminsearch}, see \cite{matlabCAT07}. To test and demonstrate the parameter estimation using the optimization Problem \ref{prob:estimationMovingBelt} we proceed as follows:
\begin{itemize}
\item Consider a process $Y_{\lambda,\sigma_1,\sigma_2}$ for given values of the parameters $\lambda,\sigma_1$ and $\sigma_2$. Compute the expected occupation times of the process corresponding to different areas $[a^i_1,b^i_1]\times[a^i_2,b^i_2]$ and belt speeds $\kappa_j$.
\item Check, if we can recover the given parameters $\lambda$, $\sigma$ and $\sigma$ by solving the above Problem \ref{prob:estimationMovingBelt}.
\end{itemize}

Let $E_{i,j}$ be the expected occupation time of the process $Y_{\lambda=1,\sigma_1=1,\sigma_2=1}$ corresponding to the time horizon $T=10$, areas $\mathcal{D}_1=[0,1.0]\times[-3.5,3.5]$, $\mathcal{D}_2=[0.5,1.5]\times[-2.5,2.5]$, $\mathcal{D}_3=[1.0,2.5]\-\times[-2.0,2.0]$ and $\mathcal{D}_4=[1.5,3.5]\-\times[-1.25\-,1.25]$ and belt speeds $\kappa_1=1$, $\kappa_2=2$. We compute the expected occupation times $E_{i,j}$  using a Monte--Carlo method based on $5000$ sample paths, see Table \ref{tab:demoMoving01}.
\begin{table} [h]
\begin{center}
\begin{tabular}{|c|c|c|c|c|}
\hline		
{\backslashbox{$\kappa_j$}{$\mathcal{D}_i$}} & $\mathcal{D}_1$	&	$\mathcal{D}_2$	&    $\mathcal{D}_3$      &       $\mathcal{D}_4$    \\  \hline
$\kappa_1=1$ 					    &	1.23698	        &      1.16451		  &	       1.63478		&  1.92876			\\ \hline
$\kappa_2=2$ 					    &	0.81322            & 	    0.70939		  &	       0.93958    &	1.06620			\\ \hline
\end{tabular}
\end{center}
\caption{Expected occupation times.}\label{tab:demoMoving01}
\end{table}
Plugging this data into Problem \ref{prob:estimationMovingBelt}, we obtain the cost functional $R(\lambda,\sigma_1,\sigma_2)$. Solving the optimization problem we recover the following values for the parameters:
\begin{equation*}
\lambda^\ast = 0.98351,\quad \sigma_1^\ast =1.01259,\quad \sigma_2^\ast =1.00878\;,
\end{equation*}
which are quite close to the originnal set $\lambda=\sigma_1=\sigma_2=1$.\\
Table \ref{tab:listOfMovingEstimations} lists for some  different settings the results of the numerical parameter estimation.

\begin{remark}
Let $0<a_1< b_1$ and $0<a_2<b_2$. Then, for a fixed belt speed $\kappa$ and time horizon $T$ we have
\begin{equation*}
{E}_{[a_1,b_1]\times[a_2,b_2],\,\kappa ,\, T}(\lambda,\sigma_1,\sigma_2)={E}_{[a_1,b_1]\times[-b_2,-a_2],\,\kappa ,\, T}(\lambda,\sigma_1,\sigma_2),
\end{equation*}
which means, the expected occupation time distribution is symmetric around the central axis of the fiber lay--down, see Figure \ref{fig:expectedOccupationTimeAreawise}. Due to this symmetry, we have just chosen symmetric domains $\mathcal{D}_i$ in the above examples. However, the selection of the areas should be adapted to the available fiber sample paths, since the occupation time distribution depends on the parameters of the lay--down process.   
\end{remark}

\begin{table}[h]
\begin{center}
\begin{tabular}{|ccc|ccc|}
\hline
\multicolumn{3}{|c|}{true parameters} & \multicolumn{3}{|c|}{recovered parameters} \\
$\lambda$ & $\sigma_1$ & $\sigma_2$ & $\lambda^\ast$ & $\sigma_1^\ast$ & $\sigma_2^\ast$    \\
\hline
0.50 & 1.00 & 1.50 & 0.50801 &  1.00003 & 1.45489\\
0.60 & 0.90 & 1.40 & 0.57202 &  0.91942 & 1.43237\\
0.70 & 2.00 & 1.75 & 0.68966 &  2.10057 & 1.77194\\
1.00 & 1.50 & 1.30 & 1.01343 & 1.48292 & 1.36154\\
1.25 & 1.25 & 2.50 & 1.23864 & 1.28534 & 2.52117\\
1.40 & 2.25 & 0.80 & 1.37193 & 2.23118 & 0.77826\\
1.50 & 0.80 & 1.75 & 1.53974 & 0.84072 & 1.73659\\
2.00 & 2.50 & 2.50 & 2.12084 & 2.52146 & 2.57379\\
2.25 & 2.50 & 1.50 & 2.11007 & 2.47998 & 1.51976\\
2.50 & 3.00 & 3.00 & 2.59092 & 3.19728 & 2.98174\\
\hline
\end{tabular}
\caption{List of recovered parameters.}\label{tab:listOfMovingEstimations}
\end{center}
\end{table}

The recovered parameters in Table \ref{tab:listOfMovingEstimations} have a maximal error of about $6\%$. This is a sufficient accuracy for many industrial applications. Nevertheless, increasing the number of sample paths used in the Monte--Carlo computations, we can improve the accuracy of the recovered parameters.

For further testing the feasibility of our estimation method we use the fiber paths simulated by FYDIST, instead of the SDE~\eqref{eqn:twoDimMovingOUmodel}. On the left hand side of Figure \ref{fig:MassDistribution}, a FYDIST fiber path is shown. Based on such fiber paths we compute the expected occupation time for different areas $\mathcal{D}_1=[0.05,0.15]\times[-0.039,0.039]$, $\mathcal{D}_2=[0.05,0.25]\times[-0.025,0.025]$, $\mathcal{D}_3=[0.1,0.20]\times[-0.018,0.018]$, $\mathcal{D}_4=[0.20,0.25]\times[-0.01,0.01]$ and $\mathcal{D}_5=[0.25,0.29]\times[-0.02,0.02]$, see Table~\ref{tab:exampleMovingFYDIST}. The belt speed $\kappa_1 = 0.0283$ and the time horizon $[0,15.93]$ are fixed.
\begin{table} [h]
\begin{center}
\begin{tabular}{|l|c|c|c|c|c|}
\hline		
{\backslashbox{$\kappa$}{$\mathcal{D}_i$}} & $\mathcal{D}_1$	&	$\mathcal{D}_2$	&    $\mathcal{D}_3$      &       $\mathcal{D}_4$ & $\mathcal{D}_5$   \\  \hline
$\kappa= 0.0283$ 					    &	5.08290	        &     7.27732		  &	      3.19642		&  1.13889		& 1.31337	\\ \hline
\end{tabular}
\caption{Expected occupation times.}\label{tab:exampleMovingFYDIST}
\end{center}
\end{table}

\begin{figure}[h]
\begin{center}
\includegraphics[ width =0.75\textwidth]{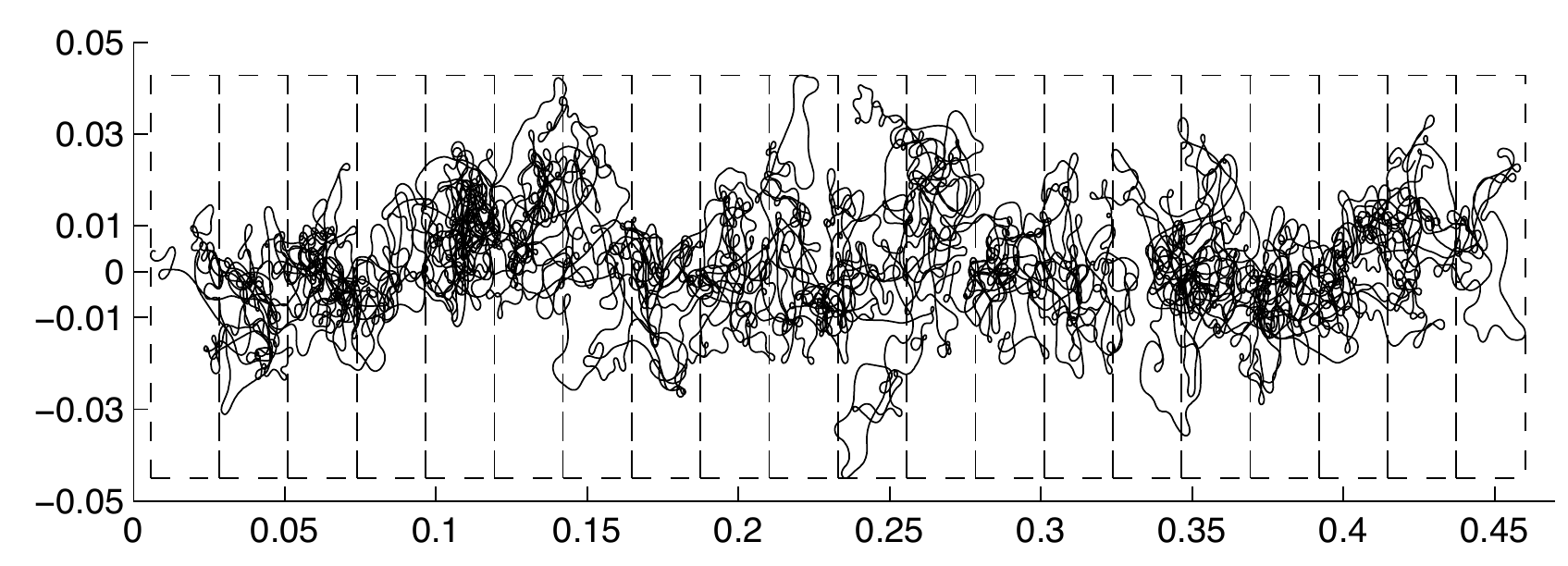}\\
\includegraphics[ width =0.75\textwidth]{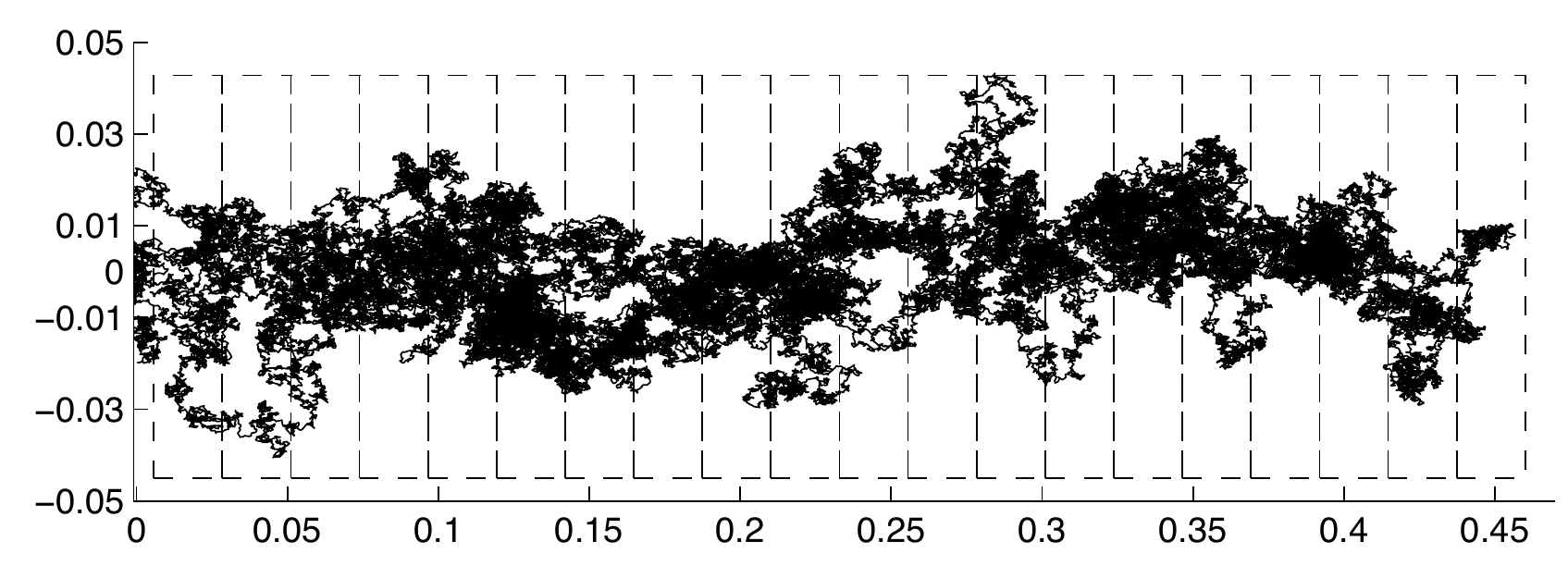}
\caption{Fiber path generated by FYDIST (top) and sample path of $Y_{\lambda^\ast,\sigma_1^\ast,\sigma_2^\ast}$ (bottom).}\label{fig:FYDISTMovingEstimatedSamples}
\includegraphics[  width =0.75\textwidth]{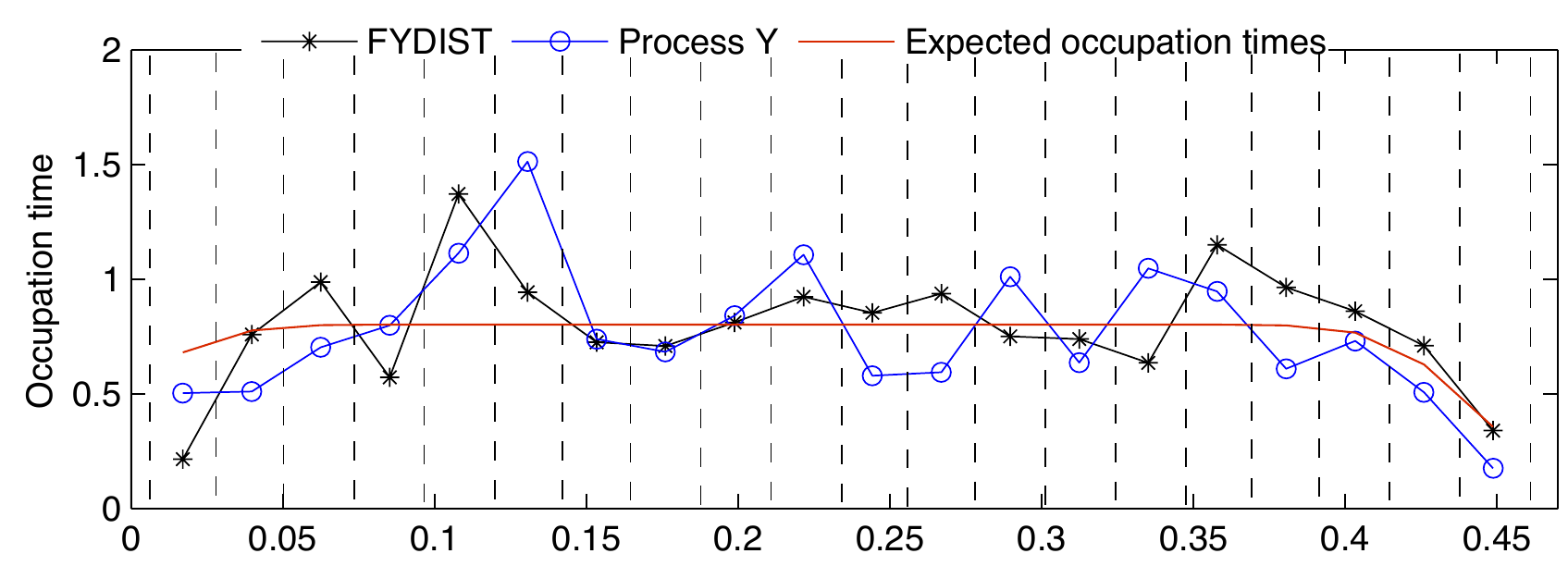}
\caption{Occupation times for the different strips and the two fiber paths shown in Figure \ref{fig:FYDISTMovingEstimatedSamples}.}\label{fig:FYDISTMovingEstimatedMassComparison}
\end{center}
\end{figure}

Solving the optimization problem, we obtain the parameters
\begin{equation}
\lambda^\ast = 4.328537,\quad\sigma_1^\ast = 0.069968, \quad\sigma_2^\ast= 0.029203\;.
\end{equation}
Figure \ref{fig:FYDISTMovingEstimatedSamples} (top) shows a FYDIST fiber path used to calculate the $E_i$ data and a sample path (bottom) of the fiber lay--down process $Y_{\lambda^\ast,\sigma_1^\ast,\sigma_2^\ast}$ characterized by the estimated parameters and the known belt speed $\kappa= 0.0283$.
%\begin{remark}
%Due to the lack of availability of FYDIST fiber samples on moving belt, in the above experiment  the expected mass distribution of FYDIST fiber paths with respect to the considered areas is not accurate. Furthermore, all the available FYDIST fiber paths are corresponded to the belt speed $\kappa=0.0283$. This lack of availability of the data decreases the accuracy of the estimated parameters.
%\end{remark}
To compare the occupation time distribution of these two fiber paths, we consider the area $[0,0.46]\times[-0.045,0.045]$ and subdivide it into $20$ vertical strips of uniform width, see Figure~\ref{fig:FYDISTMovingEstimatedSamples}. Next, we compute the occupation time of both fiber paths with respect to each strip, see Figure~\ref{fig:FYDISTMovingEstimatedMassComparison}.

\section*{Conclusion}

Based on a linear Ornstein--Uhlenbeck model~\eqref{eqn:twoDimMovingOUmodel} for the industrial fiber lay--down process, we have investigated the parameter estimation problem. For this model, we derived an analytical expression for the expected occupation time. This occupation time~\eqref{eqn:ExpOccuTimeMoving} can be regarded as the mathematical equivalent to the mass per area density of the fiber web. Given the occupation time for different domains and different belt speeds, we were able to identify the parameters of the underlying Ornstein--Uhlenbeck process by solving the related least--squares minimization problem~\eqref{eqn:costFunctionMovingBelt}.
Numerical computations based on Monte--Carlo simulations show the applicability of our method. We tested our method for fiber webs generated by an Ornstein--Uhlenbeck process as well as for webs generated using the industrial software FYDIST. Future research may focus on webs formed by more than one fiber and on nonlinear deposition models.

\bibliographystyle{abbrv}
\bibliography{biblio.bib}

\end{document}